\documentclass[12pt]{amsart}

\title[]{Differential inequalities of continuous functions and
removing  singularities  of Rado type for J-holomorphic maps }
\author[]{Xianghong Gong \and Jean-Pierre Rosay}
  \address{Department of Mathematics,
University of Wisconsin, Madison, WI 53706, U.S.A.}
\email{gong@math.wisc.edu} \email{jrosay@math.wisc.edu}

  \thanks{Research of both authors are
supported in part by NSF grants}

\newtheorem{thm}{Theorem}[section]
\newtheorem{cor}[thm]{Corollary}
\newtheorem{prop}[thm]{Proposition}
\newtheorem{lemma}[thm]{Lemma}

\theoremstyle{definition}

\newtheorem{exmp}[thm]{Example}
\newtheorem{rem}[thm]{Remark}

\newcommand{\gan}[1]{\begin{gather*} #1\end{gather*}}
\newcommand{\al}[1]{\begin{align} #1 \end{align}}
\newcommand{\aln}[1]{\begin{align*} #1 \end{align*}}
\newcommand{\eq}[2]{\begin{equation}\label{#1} #2 \end{equation}}

\renewcommand{\le}[2]{\begin{lemma}\label{lemma:#1} #2 \end{lemma}}

\newcommand{\pr}[2]{\begin{prop}\label{prop:#1} #2 \end{prop}}

\newcommand{\bp}{\begin{proof}}
\newcommand{\ep}{\end{proof}}

\newcommand{\ci}{~\cite}

\newcommand{\sk}{\medskip\noindent}

\newcommand{\myrems}{\noindent{\bf Remarks.\hspace{.5em}}}

\newcommand{\dd}[2]{\frac{\partial #1}{\partial #2}}

\newcommand{\df}{\overset{\text{def}}{=\!\!=}}
\newcommand{\D}[1]{\frac{\partial}{\partial #1}}

\newcommand{\cc}{{\bf C}}

\newcommand{\rr}{{\bf R}}

\newcommand{\rtn}{{\bf R}^{2n}}

\newcommand{\ov}{\overline}
\newcommand{\db}{\overline\partial}

\newcommand{\RE}{\operatorname{Re}}
\newcommand{\IM}{\operatorname{Im}}

\newcommand{\vol}{\operatorname{vol}}

\newcommand{\disc}{{\mathbb D}}

\newcommand{\pz}{{\,'\!z}}
\newcommand{\pw}{{\,'\!w}}

\newcommand{\dff}{\frac{\overline\partial f}{f}}

\newcommand{\fn}{function}

\newcommand{\hol}{holomorphic}
\newcommand{\jhol}{$J$-holomorphic}
\newcommand{\jholm}{$J$-holomorphic map}
\newcommand{\jholfn}{$J$-holomorphic function}
\newcommand{\jholc}{$J$-holomorphic curve}

\newcommand{\psholc}{pseudo-holomorphic curve}
\newcommand{\acs}{almost complex structure}
\newcommand{\acm}{almost complex manifold}
\newcommand{\subfn}{subharmonic  function}
\newcommand{\sub}{subharmonic}
\newcommand{\jpsh}{$J$-plurisubharmonic}
\newcommand{\jpshfn}{$J$-plurisubharmonic function}
\newcommand{\holfn}{holomorphic function}
\newcommand{\holm}{holomorphic map}
\newcommand{\holc}{holomorphic curve}
\newcommand{\contfn}{continuous function}
\newcommand{\contm}{continuous map}
\newcommand{\cont}{continuous}

\newcommand{\re}[1]{(\ref{#1})}

\newcommand{\rl}[1]{Lemma~$\ref{lemma:#1}$}

\newcommand{\sC}{\mathcal C}

\newcounter{pp}
\newcommand{\bpp}{\begin{list}{$\hspace{-1em}\arabic{pp}.$}{\usecounter{pp}}}
\newcommand{\epp}{\end{list}}

\newcounter{ppp}
\newcommand{\bppp}{\begin{list}{$\hspace{-1em}(\roman{ppp})$}{\usecounter{ppp}}}
\newcommand{\eppp}{\end{list}}

\begin{document}
\maketitle
\setcounter{thm}{0}\setcounter{equation}{0}

\centerline{\bf Introduction}

\

Among the results of Part One  we have:

\newcommand{\n}{$\ \ \clubsuit\clubsuit\vartriangleright\ \ $}
\newcommand{\nb}{$\ \ \vartriangleleft\clubsuit\clubsuit\ \ $}

\sk {\bf Proposition A}. {\em Let $\Omega$ be an open set in
$\cc^n$ and $f\colon\Omega\to\cc$ be a \contm.   If   on
$\Omega\setminus f^{-1}(0)$ $|\db f|\leq K|f|$ ($K$ a positive
constant), then $f^{-1}(0)$ is an analytic set.}

\medskip

Among the results of Part Two there is the following:

\sk  {\bf Proposition B}. {\em Let $J$ be a $C^1$-smooth \acs\
defined
  in $\rtn$. Let  $v$ be a \contm\ from the unit disc $\disc\subset
  \cc$ into $(\rtn,J)$. Let $u$ be either a constant map
  from $\disc$ into $\rtn$
  or a proper
  \jhol\ and $\mathcal C^2$-smooth  map from $\disc$ into an open subset
$\Omega$ of
   $\rtn$.
   Assume that $v$ is \jhol\ near $z$ if $v(z)\not\in u(\disc)$.
  If $u\equiv u(0)$ then $v$ is \jhol\ on $\disc$;
  if  $u\not\equiv u(0)$ and $v(0)\in u(\disc)$,
   either $v$ maps a neighborhood of $0$ into $u(\disc)$, or $v$ is
  \jhol\ near $0$.}

  If $u\equiv u(0)$  the hypothesis is thereof that $v$ is \jhol\ at any point
  $z$ such that $v(z)\neq u(0)$. Proposition B generalizes the classical Rado
  Theorem.

\medskip
The more general statements in the text that contain
Propositions~A and B are Theorem~A and  Theorem~C.  

 Parts  One and Two are independent. A unifying theme is the
elimination of exceptional sets on which no a-priori assumption is
to be made, in the style of Rado's Theorem. In some of our
statements    this   
 is the only novelty.
  Another common feature is the
role played by differential inequality of the type $|\db f|\leq
|f|$ or $|\db f| \leq\epsilon|\partial  f|$.

\bigskip
\begin{center}
Organization of the paper\\
- Part One -
\end{center}

\sk
1) $|\db f|\leq |f|$ versus $\db f=Af$. Analyticity of zero set
of $f$
  and   a  counter-example

\sk 2)    One dimensional vector-valued case

\sk 3) Proof  of  
 the analyticity of  $f^{-1}(0)$
(scalar-valued $f$ defined on $\cc^n$)

\sk
4) Remarks on uniqueness results on \jholc s

\newpage
\sk
\begin{center}
- Part Two -
\end{center}

\sk
1) Removability of polar sets for \jholc s

\sk
2) Rado type theorem for \jholc s

\

\setcounter{thm}{0}\setcounter{equation}{0}

\begin{center}
\bf Part One: Analyticity of zero sets of  \contfn
\end{center}

\newcommand{\ost}{\Omega\setminus f^{-1}(0)}

\sk
{\bf 1. $\db$-differential inequalities and  results of Part One}

\

Let $f$ be a \contfn\ on an open subset $\Omega$ of $\cc$. Denote $\Omega^*
=\Omega\setminus f^{-1}(0)$.
By the $\db$-differential inequality $|\db f|\leq K|f|$  on
$\Omega^*$, we mean that on $\Omega^*$,
$\db f=Af$ holds
in the distribution sense, and $A\in L^\infty_{(0,1)}(\Omega^*)$.

In general, one may,  of course,
  consider a \contm\ $f\colon \Omega\to\cc^m$ satisfying
$\db f=Af$ on $\Omega^*=\Omega\setminus f^{-1}(0)$, where $\Omega$
is now an open subset of $\cc^n$ and
$A$ is an $m\times m$ matrix of $(0,1)$-forms on   $\Omega^*$. We are
interested in
the analyticity of the zero set $f^{-1}(0)$.

\sk

We formulate our results in following three cases: i) $f$
is a vector-valued \cont\ function   on $\cc$, ii) $f$ is a \cont\ function on
$\cc^n$,
iii) $f$
is a vector-valued \cont\ function on $\cc^n$ under the additional assumption
that $f^{-1}(0)$ is already real analytic.

\sk

The result for case (iii) will be given later by Proposition~3.4.  Here we describe
results for (i) and (ii).

\

i) Vector-valued $f$, defined on $\cc$.

\sk

Our first result  is the following:

\noindent {\bf Proposition C.} {\em Let $\Omega$ be   a
connected open subset of $\cc$, and let
   $f\colon\Omega\to\cc^m$ be a \contm. Assume that $\db f=Af$ on $\ost$,   where
$A$ is an $m\times m$ matrix with entries in  $
L^p_{(0,1)}(\ost)$. If $2<p\leq\infty$ and $f\not\equiv0$,
  then $f$ has isolated zeros in $\Omega$ and $\db f= Af$ holds on $\Omega$,
  where   $A$ is set to be zero on $f^{-1}(0)$.
}

\

ii) Scalar-valued $f$, defined on $\cc^n$.

\sk

Next we turn to \contfn s $f$  on an open subset $\Omega$ of
$\cc^n$ with $n>1$. Our approach is based on the above
one-dimensional result, by applying it to all  complex lines $L$
parallel to a coordinate axis  in $\cc^n$. Thus one needs to
understand the tedious question about the restriction of $\db
f=Af$ to all $z_j$-lines, including the restriction of $|\db
f|\leq K|f|$ to lines, needed for Proposition A.   (Of course,
$\db (f|L)$ makes sense since $f$ is \cont.)

  We
postpone
  this technicality to a later discussion.
  Here we ask that for every complex line $L$
   (intersecting $\Omega$ and)  parallel to a coordinate axis,
   in the distribution sense
$$
\db (f|_L)=A_Lf,   \ \text{on $L\cap \Omega\setminus f^{-1}(0)$}, \quad
  \|A_L\|_{L^p(L\cap\Omega\setminus f^{-1}(0))}\leq M,
 $$
where $M>0$ is a constant   (independent of $L$).

\

Equivalently, after we prove Proposition C  in section 2, the
above equation can be assumed to hold on the line (not only off the
zero set),   by setting $A|_L=0$ on $f^{-1}(0)$ and by assuming
$p>2$. Therefore  we restate it in a simpler form:
$$
\db (f|_L)=A_Lf   \ \text{on $L\cap \Omega$}, \quad
  \|A_L\|_{L^p(L\cap\Omega)}\leq M
 \leqno{(L)} $$
for a constant $M$ independent of $L$.

\

We  now state the main result of Part One:


\sk
{\bf Theorem A.} {\em Let $2<p\leq\infty$.
Let $f$ be a \contfn\ on an  open subset
$\Omega$ of $\cc^n$. Assume that   for each $j=1,\dots, n$,
 $f$ satisfies (L) for all complex lines   $L$   parallel to
the $z_j$-axis.
Then $f^{-1}(0)$ is a complex variety.
Moreover, $e^{-u}f$ is \hol, for some
  locally defined
$u\in C_{loc}^{\alpha}$,
with $\alpha=1-\frac{2}{p}$ for $2<p<\infty$, and any $\alpha<1$ for
$p=\infty$.
}

\medskip

Condition (L) plays two roles in the crucial \rl{n1}: (a) it is
used to show that the trivial extension of $A$ is $\db$-closed,
(b) once this is proved it is used again for proving that the
solution to some $\db$-problem is $C_{loc}^\alpha$-smooth.

To show the analyticity of $f^{-1}(0)$, it suffices to find $u\in
L_{loc}^{\infty}(\Omega)$ solving $\db u=\dff$ on $\Omega\setminus
f^{-1}(0)$, or on $\Omega$ by extending $\dff$ trivially by $0$ on
$f^{-1}(0)$. There are other ways to insure boundedness, or
continuity, or   $C_{loc}^\alpha$-smoothness   of   solutions to
the $\db$-problem. This is discussed in details at the beginning
of Section 3.

\bigskip\noindent
{\bf Discussion and an open problem.}

\sk

The next example shows that the above two results do not hold for $p\leq2$.

\sk
{\bf Example.}
Let $a_k\in\cc^n$ with
$0<|a_k|<1/2$ and $a_k\to 0$ as $k\to\infty$. Let
$$
f(z)=\frac{1}{(\log |z|)\prod_{k=1}^\infty \bigl|\log
|z-a_k|\bigr|^{\frac{1}{k^2}}}.
$$
Then $f$ is \cont\ on $\disc_{1/2}^n$, and $\db (f_L)=a_Lf$ with
$$
a_L(z)=-\frac{z\cdot d\ov z}{2|z|^2\log |z|}-\sum\frac{1}{2k^2}
\frac{(z-a_k)\cdot d\ov z}{|z-a_k|^2\log |z-a_k|}\in L^2(L\cap\disc_{1/2}^n).
$$
Also $\db f=Af$ with
$A\in L^{2n}(\disc_{1/2}^n)$.

\sk

Theorem A is motivated by recent work of  Pali\ci{Pazefi},
who proved, among other results,
that  $f^{-1}(0)$ is a complex variety
under the assumptions that $f$ is a vector of  smooth \fn s on $\Omega\subset\cc^n$,
$\db f=Af$ on $\Omega$ with $A\in C^\infty$, and $f$ satisfies a certain finite
resolution property.
In\ci{Pazefi}, the difficulties lie in the case   when
   $f$ is not a scalar function.

\sk

The following problem remains open:

\sk
{\bf Problem A.} (1) Let $f$ be a  vector of \contfn s on a domain
$\Omega\subset\cc^n$ satisfying the $\db$-inequality
  $|\db f|\leq K|f|$
on $\Omega\setminus f^{-1}(0)$. Is $f^{-1}(0)$ a complex
variety?


(2) Back to the scalar case: Let $f$ be a \contfn\ on
$\Omega\subset\cc^n (n>1)$. Assume that $\db f = Af$ on
$\Omega\setminus f^{-1}(0)$ with $A\in L_{(0,1)}^p(\Omega\setminus
f^{-1}(0))$. If $2n<p<\infty$, is $f^{-1}(0)$ a complex variety?

\sk

\setcounter{thm}{0}\setcounter{equation}{1}
\setcounter{section}{2}

\newpage

\noindent {\bf 2.  One dimensional vector-valued case}

\

Let $\disc$ be the unit disc in $\cc$,
   $\disc_r$ the disc of radius $r$ centered at the origin, and
   $\disc_r(z)$   the disc  with  radius $r$,
centered at $z$. For $f\in L^p(\disc_\delta)$ with
$2<p\leq\infty$, define
$$
Tf(z)=\frac{i}{2\pi}\int_{|\xi|<\delta}\frac{f(\xi)}{z-\xi}
 d\xi\wedge d\ov\xi.
$$
Recall the well-known estimates
 \gan{
|Tf(z)|\leq c_p\delta^{1-\frac{2}{p}}\|f\|_{L^p(\disc_\delta)},\quad 2< p\leq\infty,\\
|Tf(z')-Tf(z)|\leq c_p\|f\|_{L^p(\disc_\delta)}|z'-z|^{1-\frac{2}{p}},\  2<p<\infty,\\
|Tf(z')-Tf(z)|\leq
c_p\|f\|_{L^\infty(\disc_\delta)}|z'-z|(1+\bigl|\log\delta\bigr|+\bigl|\log|z'-z|\bigr|),
}
where $c_p$ depends only on $p$.
Also $\dd{}{\ov z}Tf=f$ on $\disc_\delta$ in the distribution sense
(see\ci{Ahsisi}).

  The following lemma for the scalar case will be repeatedly used.
Note that the lemma gives us
  Proposition~C.
\le{1}{  Let $2<p\leq\infty$.  Let
$\phi_p(\delta)=\delta^{1-\frac{2}{p}}$ for $p\neq\infty$ and
$\phi_\infty(\delta)=\delta(1+|\log\delta|)$.   Let  $f\colon
\disc\to \cc^m$ be a \contm.
  Assume that for some $m\times m$
matrix $A=(a_{jk})$ with $a_{jk}\in L^p(\disc\setminus f^{-1}(0))$
$$
\dd{f(z)}{\ov z}= A(z)f(z), \quad \|A\|_{L^p(\disc\setminus f^{-1}(0))}\leq
M<\infty
$$
holds
in the distribution sense
on $\disc\setminus f^{-1}(0)$. There is a constant $\epsilon>0$ such that if
$\delta>0$ and
    $M\phi_p(\delta) <\epsilon$,   then for each
$z\in\disc$, there is an invertible $m\times m$ matrix $I+g$
  on $\disc_{\delta}(z)$
  such that
$(I+g)f$ is holomorphic on $\disc_{\delta}(z)\cap\disc$,   and
$g\in C^\alpha(\disc_\delta(z))$ with  $\alpha=1-\frac{2}{p}$ for
$2<p<\infty$ and any $\alpha<1$ for $p=\infty$. Moreover,
$\sup_{w\in\disc_\delta(z)}|g(w)| \leq c_{p,m}M\phi_p(\delta)$
 for some constant   $c_{p,m}$,  and
$\dd{f}{\ov z}=\tilde Af$
on $\disc$ where $\tilde A$ equals $A$ on $\disc\setminus f^{-1}(0)$ and  zero on $f^{-1}(0)$.
}
\begin{proof} Set $\tilde A=A$ on $\disc\setminus f^{-1}(0)$
and $\tilde A=0$ on $f^{-1}(0)\cup (\cc\setminus \disc)$.
Still denote $\tilde A$ by $A$.
We may assume that $z=0$.
  We need
$$
0=\D{\ov z}((I+g)f)=(\D{\ov z}g+(I+g)A)f
$$
on $\disc_\delta$.  For a matrix $u$ of
  functions in $L^p(\disc_\delta)$,   define
$$
Tu(z)= \frac{i}{2\pi}\int_{|\xi|<\delta}\frac{u(\xi)}{z-\xi}\, d\xi\wedge d\ov \xi,
$$
where the integration takes place entry by entry. The $\delta$,
depending only on $m$ and $p$,
will be determined later.

Consider the equation
$$
   g+T((I+g)A)=0.
\leqno{(\ast)}
$$
  For the matrices $g$   and $A=(a_{jk})$, put
$$
\|g\|=\sup_{j,k,\disc_\delta}\{|g_{jk}(z)|\}, \quad
\|A\|_p=\max_{j,k}\{\|a_{jk}\|_{L^p(\disc_\delta)}\}.
$$
Recall that $T((I+g)A)$ is \cont\ on $\ov\disc_\delta$, if $g$ is \cont\ on
$\ov\disc_\delta$.
We also have
\gan{
\|T((I+g_2)A)-T((I+g_1)A)\|\leq  c_{p,m}\phi_p(\delta)\|A\|_p\cdot\|g_2-g_1\|,\\
\|T(A)\|\leq c_{p,m}\phi_p(\delta)\|A\|_p. } Therefore,
$g=\lim_{k\to\infty}S^k(0)$, where $S(g)=-T((I+g)A)$,
 is the unique matrix
   that is \cont\ on $\ov\disc_\delta$
     and satisfies  ($\ast$) on $\disc_\delta$,
     provided
      $$
M\phi_p(\delta)<\frac{1}{2c_{p,m}}.
$$
      Moreover, the $g$
    satisfies
$$
\|g\|\leq \frac{c_{p,m}\|A\|_p\phi_p(\delta)}{1-c_{p,m}\|A\|_p\phi_p(\delta)}\leq
2c_{p,m}M\phi_p(\delta).
$$
By ($\ast$), $(I+g)f$ is \hol\ on $\disc_\delta\setminus f^{-1}(0)$.
Let $h$ be an entry of
$(I+g)f$.
We know that $h$ is \cont\ on $\disc_\delta$. If $h$ does not vanish at
$z_0\in\disc_\delta$, then $f(z_0)\neq0$ either. Hence $h$ is \hol\ near $z_0$.
By
Rado's theorem, $h$ is holomorphic on $ \disc_\delta$.
  Since $g=-T((I+g)A)$ and $(I+g)A$ is in $L_{loc}^p$ then $g\in C_{loc}^\alpha$
for $\alpha=1-\frac{2}{p}$ when $2<p<\infty$ and any $\alpha<1$ for $p=\infty$.

Since $f$ is identically zero or has isolated zeros, it is easy to see that
$\dd{f}{\ov z}=\tilde Af$ on $\disc$.
\end{proof}

  When $f$ is scalar, one simply takes the above $1+g$  to be the so-called
 {\em integrating factor}, $e^{-u}$, where
$u=A\ast\frac{1}{\pi z}$ by setting  $A=0$ on
$f^{-1}(0)$. On the other hand,  the   above   contraction
argument avoids the Rado theorem, if one has a stronger
condition that $\dd{ f}{\ov z}=Af$ on $D$.

The zero set $f^{-1}(0)$ for the one-dimensional case is
  studied by Ivashkovich-Shevchishin\ci{ISnini}, under the assumption
  that $\db f=Af$ on the whole domain $\Omega\subset\cc$. See also
  Floer-Hofer-Salamon\ci{FHSnifi}.


\setcounter{thm}{0}\setcounter{equation}{0}
\setcounter{section}{3}

\

\noindent
{\bf 3. Proof of Theorem A (scalar-valued case, functions defined on $\cc^n$)}

\

The heart of the matter is to prove that $B$, the trivial extension of
${\overline\partial f \over f}$, is $\overline\partial$-closed.
This is done in Lemmas 3.2 and 3.3
but we start in Lemma~3.1  with some remarks that will allow us later to
prove that the solutions to $\db u=B$ are locally  bounded.

\sk
{\bf 3.1. Restriction of $\db f=Af$ on lines.}

\sk

First we need to discuss
  the restriction of $\db f=Af$ on complex  lines.

However, the reader could skip   this section 3.1 (remarks a)-d) below) and go directly to section 3.2
(beginning of the
proof of Theorem A),   by
first accepting
    the assertion:  if   $f$ satisfies condition (L) for almost all $z_j$-lines
   and if $B$, the trivial extension
of $\dff$ by $0$ on $f^{-1}(0)$, is $\db$-closed,  then all solutions   $u$ to  $\db u=B$
are locally bounded on $\Omega$. The latter will suffice for  the proof of analyticity
of $f^{-1}(0)$.

\sk

Let $f$ be a \contfn\ defined on an open subset
$D$ of $\cc^n$. It is worth pointing out the following.

a) Let $1<p\leq\infty$.
  Let $M$ be a positive constant. Assume that for almost all $z_1$-lines $L$
$$
\dd{f|_L}{\ov z_1}=A_Lf, \qquad \|A_L\|_{L^p(L\cap D)}\leq
M.\leqno{(L)}
$$
Then in the distribution sense
$$
\dd{f}{\ov z_1}=Bf
$$
on $D$, for some function
$B\in L^p_{loc}(D)$, that we set to be 0 on $f^{-1}(0)$.

That function $B$ satisfies:
$$  \| B\|_{L^p(D\cap L_\epsilon)}\leq
(\pi\epsilon^2)^{n-1\over p}M,\leqno{(E)}$$
where  $L_\epsilon$ is the $\epsilon$-neighborhood of
$L$
in $\cc^n$ equipped with the sup-norm. (In the latter formulation  we avoid
the problem of having to restrict $B$ to lines.)

\

\noindent{\em Proof.} Check $\dd{f}{\ov z_1}=B f$ for some
$B\in L^p_{loc}(D)$.   Fix a compact subset $K$ of $D$. Let
$\varphi$ be a smooth function whose support is contained  in  $K$.
Put $z'=(z_2,\dots,z_n)$ and
$$
dV(z')=(\frac{i}{2})^{n-1}dz_2\wedge d\ov z_2\wedge
\cdots\wedge dz_{n}\wedge d\ov z_n, \quad
dV(z)=\frac{i}{2}dz_1\wedge d\ov z_1\wedge
dV(z').
$$
Then \aln{ |(\dd{f}{\ov z_1},\varphi)|&=|\int\int
f\dd{\varphi}{\ov z_1}
\, \frac{i}{2}dz_1\wedge d\ov z_1\wedge dV(z')|\\
&=|\int( \dd{f|_{\text{$z'$ fixed}}}{\ov z_1},\varphi)_{z_1}\, dV(z')|\\
&=|\int(\int A|_{L(z')}f\varphi\, \frac{i}{2}dz_1\wedge d\ov z_1)
 dV(z')|, \ L(z')=\{(z_1,
z')\colon z_1\in\cc\}\cap D\\
&\underset{\text{H\"older}}{\leq}M\int\|f\varphi\|_{L^q(z_1)}\, dV(z')
\underset{\text{H\"older}}{\leq} M|K'|^{\frac{1}{p}}\left[\int |f\varphi|^q
\frac{i}{2}dz_1\wedge d\ov z_1\wedge dV(z')\right]^{\frac{1}{q}},
}
where $K'$ is the projection of $K$ by $z\to(z_2,\dots, z_n)$   and
$|K'|=\int_{K'}dV(z')$.
So $|(\dd{f}{\ov z_1},\varphi)|\leq M|K'|^{\frac{1}{p}}\|f\varphi\|_{L^q}$, i.e
$f\varphi
\to(\dd{f}{\ov z_1},\varphi)$ is \cont\ in $L^q$-norm.
So $(\dd{f}{\ov z_1},\varphi)=\int B f\varphi\, dV(z)$
for some $B$ supported in $K$,   vanishing on $f^{-1}(0)$,   and satisfying
\eq{alpl+}
{\|B\|_{L^p}\leq M|K'|^{\frac{1}{p}}.
}
For the general case, let $\{\chi_k\}$ be a partition of unity for a locally
finite open covering $\{U_j\}$ of $D$. Then $\dd{f}{\ov z_1}=B_kf$ on $U_k$,
and   $(\dd{f}{\ov z_1},\varphi)=\sum(\dd{f}{\ov
z_1},\chi_k\varphi)=(B,\varphi)$
for $B=\sum\chi_kB_k\in L^p_{loc}$.   Finally, (E)  comes from \re{alpl+}
and the uniqueness of $B$ (up to a set of measure $0$ and   requiring
$B|_{f^{-1}(0)}=0$) since $f$
is \cont.
\hfill {QED}

\sk

b) Let $1<p\leq\infty$. Assume that  $\dd{f}{\ov z_1}=Bf$ on $D$  with $B$
satisfying (E)
for all $z_1$-lines $L$ and all $\epsilon>0$.
 Then ($\text{L}$) holds for all $z_1$-lines $L$.
In other words,  when $f$ is  \cont,  ($\text{L}$) and (E)  are
equivalent. Moreover,  that ($\text{L}$) holds    for almost
all $z_1$-lines implies ($\text{L}$)   for all $z_1$-lines.
\bp Let $\chi\geq0$ be a smooth function   on $\cc^n$ with compact
support and   $\int \chi(z) dV(z)=1$. Put
$\chi_\delta(z)=\frac{1}{\delta^{2n}}\chi(\delta z)$ and
$B_\delta=B\ast\chi_\delta$. We may assume that $L$ is
  the $z_1$-axis.  Fix a compact subset $K$ of
   $L\equiv L\cap D\subset\cc$.
Then
\al{\label{ev}
\|B_\delta\|_{L^{p}(K\times \disc_\epsilon^{n-1})}&=\left(\int_{K\times\disc_\epsilon^{n-1}}|\int B(z-\delta w)\chi(w)\, dV(w)|^p\,
dV(z)\right)^{\frac{1}{p}}\\
&\leq\int\left(\int_{K\times\disc_\epsilon^{n-1}}|B(z-\delta w)|^p\,
dV(z)\right)^{\frac{1}{p}}\chi(w)\,  dV(w)
\nonumber\\
&\leq\int\left(\int_{D\cap\cc\times(-\delta w'+\disc_\epsilon^{n-1})}|B(z
)|^p\, dV(z)\right)^{\frac{1}{p}}\chi(w)\,
dV(w)
\nonumber\\
&\leq
(\pi\epsilon^2)^{\frac{n-1}{p}}M, \qquad w=(w_1,w').\nonumber
}
Thus, $\|B_\delta|_{K\times\{0'\}}\|_{L^p}\leq M$, since $B_\delta$ is \cont.
For each smooth function
$\varphi$ with compact support in $L$,
$$(\dd{f|_{L}}{\ov z_1},\varphi)=\lim_{\delta\to0}
(\dd{{f_\delta}|_{L}}{\ov z_1},\varphi)=\lim_{\delta\to0}((Bf)_\delta|_L,\varphi)=\lim_{\delta\to0}(B_\delta|_Lf,\varphi).
$$
  For small $\delta>0$
$$\|B_\delta(\cdot,z')\|_{L^p(K)}=\left(\int_K|\int B((z_1,z')-\delta
w)\chi(w)\, dV(w)|^p\, \frac{i}{2}dz_1\wedge d\ov z_1\right)^{\frac{1}{p}}
\leq M.$$ Using
the weak compactness in $L^p$, we find  a sequence $B_{\delta_k}|_L$
  that   converges weakly to some $A_L\in L^p(L)$ with $\|A_L\|_{L^p}\leq M$.
\ep

 \sk

  By Proposition C, remarks a) and b), we know that the condition
$$
\db f=Af \quad \text{on $\Omega\setminus f^{-1}(0)$}, \quad A\in L^\infty(\Omega\setminus f^{-1}(0))
$$
is equivalent to
$$
\db f=Af \quad \text{on $\Omega$}, \quad A\in L^\infty(\Omega),
$$
 in which $A$ is set to be $0$ on $f^{-1}(0)$.

 \

c)   Let $v=v_1\,d\ov z_1+\cdots+v_n\, d\ov z_n$ with
$v_j\in L^1_{loc}(\Omega)$ be  $\db$-closed   on $D$. The local existence and  boundedness
  of solutions
$u$ to $\db u=v$ follow
from any of the following
conditions:

\smallskip
(i) The integral of $\int_D\frac{|v(w)|}{|z-w|^{2n-1}}\, dV(w)$ is bounded (independently of $z$).

\smallskip (ii) For each
$v_j$ there is $k$ such that  for some $p>2$
and $M>0$, $\|v_j\|_{L^p(D\cap L)}\leq M$  for
(almost) all $z_k$ -lines $L$.

\smallskip(iii)  $v$ is in $L^p(D)$ for some $p>2n$.

\smallskip
In fact,  all local solutions are \cont\ for case (ii), and are in $C_{loc}^\alpha$
with $\alpha=1-\frac{2n}{p}$ for $2n<p<\infty$ and any $\alpha<1$ for $p=\infty$, for
(iii).
\bp We may assume that   $D=\disc^n$.
It is elementary that (iii) implies (i). (ii) implies (i) too:
 Assume that $\|v_j\|_{L^p(D\cap L)}\leq M$
for all $z_k$-lines $L$ (say $k=n$).  Then
\aln{
&\int_D\frac{|v_j(w)|}{|z-w|^{2n-1}}\, dV(w)
\leq M\int_{|\pz|<1}\left(\int_{|z_n|<1}
\frac{1}{(|\pz|^2+|z_n|^2)^{(2n-1)q/2}}\, \frac{i}{2}dz_n\wedge d\ov z_n
\right)^{\frac{1}{q}}\, dV(\pz)\\  &\qquad
\leq cM\int_{|\pz|<1}1+|\pz|^{\frac{2}{q}-(2n-1)}\, dV(\pz)=\tilde c M\int_0^1
r^{2n-3}+r^{-\frac{2}{p}}\, dr<\infty
}
for $\frac{1}{q}=1-\frac{1}{p}$ (and $n>1$).

Assume that (i) holds with $D$ being the unit ball
and $\db v=0$.

Recall the Bochner-Martinelli   kernels:
$$
\frac{(n-1)!}{(2\pi i)^n}\sum_j(-1)^{j-1}\frac{\ov\zeta_j-\ov
z_j}{|\zeta-z|^{2n}}\wedge_{k\neq j}(d\ov \zeta_k-d\ov z_k)
  \wedge d \zeta_1\wedge\cdots\wedge d\zeta_n=\sum_{q=0}^{n-1}\Omega_q(z,\zeta),
$$
where $\Omega_q$ is  of bidegree $(n,n-q-1)$ in $\zeta$ and
 bidegree $(0,q)$ in $z$. For a
$(0,q)$-form $u$ on   $D$ with $q>0$, define
$$
 B_M\ast u(z)=-\int_{\zeta\in D}u(\zeta)\wedge \Omega_{q-1}(z,\zeta).
$$
By the Koppelman formula ([H-L], p. 57), $w=\db (B_M\ast w)-B_M\ast\db w$, if
$w$ is a smooth $(0,1)$-form on $D$ with compact support.

Let $v_\epsilon=\sum_{j=1}^n(v_j)_\epsilon d\ov z_j$,
where the smoothing $(v_j)_\epsilon$ is defined in the proof of remark b).
Let $\chi$ be a smooth function supported in the ball $D_{1/2}\colon |z|<1/2$.
For small and positive $\epsilon$, $v_\epsilon$ is well-defined and $\db$-closed
on $D_{1/2}$. Thus
$$
 \chi v_\epsilon=\db (B_M \ast (\chi v_\epsilon))-B_M\ast (\db\chi \wedge v_\epsilon),
  \quad \text{on $D$}.
$$
By the Fubini theorem, $\|B_M\ast(\chi v_\epsilon)\|_{L^\infty(D_{1/2})}+
\|B_M\ast (\db\chi \wedge v_\epsilon)\|_{L^\infty(D_{1/2})}\leq c
\||v|\ast\frac{1}{|z|^{2n-1}}\|_{L^\infty(D)}$ for some constant $c$.
On $\ov D_{1/3}$, $B_M\ast (\db\chi \wedge v_\epsilon)$ is
smooth and $\db$-closed for small $\epsilon>0$. It is well-known that there is a smooth function $\tilde u_\epsilon$
such that $\db \tilde u_\epsilon=B_M\ast (\db\chi \wedge v_\epsilon)$
and $\|\tilde u_\epsilon\|_{L^\infty(D_{1/3})}\leq \tilde c\|B_M\ast
(\db\chi \wedge v_\epsilon)\|_{L^\infty(D_{1/3})}$.
We now have $\db (B_M \ast (\chi v_\epsilon)-\tilde u_\epsilon)=v_\epsilon$ on $D_{1/3}$.
By the weak compactness, the week-limit $u$ of some sequence
$B_M \ast (\chi v_{\epsilon_k})-\tilde u_{\epsilon_k}$ is  bounded and satisfies $\db u=v$.

For (ii),  the continuity of $u$ follows from the continuity  of
$|v_j|\ast\frac{1}{|z|^{2n-1}}$. For (iii), it is elementary that
$B_M\ast v\in C_{loc}^\alpha$
when $v\in L^p(D)$ and $p>2n$.
\ep

For a complex line $L$ in $\cc^n$, denote by $L_\epsilon$ the $\epsilon$-neighborhood
of $L$
in $\cc^n$    with the  sup-norm.

\sk

We now point out a special case, which is actually used in our proofs.

\medskip d) Let $2<p\leq\infty$.
 Let $u\in L^1_{loc}(D)$ and $\db u=v_1\, d\ov z_1+\cdots+v_n\, d\ov z_n$. Assume that
for each $j$ there is a constant $M$ such that
$$
\|v_j\|_{L^p(D\cap L_{j,\epsilon})}\leq M(\pi\epsilon^2)^{\frac{n-1}{p}}
$$
for  all $z_j$-lines $L_j$ and all $\epsilon>0$,
then $u\in C^\alpha_{loc}$ for $\alpha=1-\frac{2}{p}$
 ($2<p<\infty$) or for all $\alpha<1$ ($p=\infty$).
\bp
By c (ii), we know that all local solutions $u$ are \cont.
Repeat the argument in b) for $\dd{u}{\ov z_j}=v_j$ (in a simpler way).
We have $\dd{u|_L}{\ov z_j}=\tilde v_L$
with $\|\tilde v_L\|_{L^p}\leq M$. Now the local H\"older-$\alpha$ norm of $u|_L$ is bounded by some constant
independent of $L$. Therefore   $u|_L$ (and hence $u$)
is $C^\alpha_{loc}$-smooth,
by the well-known
estimates
stated at the beginning of Section 2.
\ep

\sk
{\bf 3.2. $\db$-closedness of trivial extension of $\dff$, when $f^{-1}(0)$ is the graph of a \contfn.}

As already said, the heart of the matter is to prove that $B$, the trivial extension of
${\overline\partial f \over f}$, is $\overline\partial$-closed.
Once this is done
the existence of an integrating factor and hence the analyticity of $f^{-1}(0)$
is trivial. However the proof goes as follows. In this section 3.2 we
prove
$\overline\partial$-closedness of
$B$ near specific points of $f^{-1}(0)$ (\rl{n1}). We then get the
analyticity of $f^{-1}(0)$ off an exceptional set (in the notations below, the set
where $N(p)=\infty$). Using induction on dimension, the proof of
$\overline\partial$-closedness of $B$ is then achieved in section~3.3.

By our assumptions,  if $f$ vanishes at $p$, then  for any complex
 line  $L$ passing through $p$ and parallel to a coordinate axis, there is a \contfn\ $u$ so that
$e^{-u}f|_L$ is holomorphic on $L\cap\Omega$.  Define the order of
vanishing  of  $f|_L$   at $p$
to be the vanishing order
of
$e^{-u}f|_L$ at $p$.
We also define the number of zeros of  $f|_L$  to be the  number of zeros of
$e^{-u}f|_L$. Both are  independent of the choice of $u$.

We start with the following.
\le{2}{Let $2<p\leq\infty$.  Let $f(z)$ be a \cont\ function on $\disc^n$.
  Assume  that there exists a positive constant $M$
such that
for all complex lines $L$ parallel to the $z_n$-axis,  $\db (f|_{L})=A_Lf$
on $\disc^n\cap L$
and
$\|A_L\|_{L^p(L\cap\disc^n)}\leq M$.
Assume that $f(0)=0$ and $f(0,z_n)$ does
not vanish when $0<|z_n|\leq\epsilon_0<1$.

(i) There exists $\delta>0$ such that for
each $\pz\in\disc_\delta^{n-1}$, $f(\pz,z_n)$ has exactly $N$
zeros $z_n= r_j(\pz), j=1,\dots, N$ in the disc $\{z_n\colon
|z_n|<\epsilon_0\}$,  counting multiplicity. Moreover,  for each
$\epsilon\in(0,\epsilon_0)$ if the $\delta$ is sufficiently small we have
$|r_j(\pz)|<\epsilon$  for $|\pz|<\delta$ and $\lim_{\pz\to0}r_j(\pz)=0$ for all
$j$.

(ii)
If the above $r_j$ are all the same, denoted by $r$, then
$f(z)=((z_n-r(\pz))v(z))^N$, where  $v$
is \cont\ on  $\disc_\delta^{n-1}\times\disc_{\epsilon_0}
\setminus f^{-1}(0)$, and $v,\frac{1}{v}$ are
locally bounded on
$\disc_\delta^{n-1}\times\disc_{\epsilon_0}$.   }
\begin{proof}
(i)
By Proposition C, we know that there exists $0<\epsilon_1\leq\epsilon_0$
such that for each $\epsilon\in(0,\epsilon_1)$ there exists a
function $g(z)$
such that $(1+g(z))f(z)$ is holomorphic in $z_n$ and
$|g(z)|\leq 1/4$ for $z\in\disc^{n-1}\times\disc_\epsilon$.
  By definition   the number of zeroes
of $f$ as a function of $z_n$ is the same as the number of zeroes of the
holomorphic function $(1+g)f$. So it is the winding number of
 $(1+g)f$ and    hence that   of $f$, which is
constant (in a neighborhood of 0).
The last conclusion follows from  the continuity of $f$ also.

(ii) We reduce it to case $N=1$ first.
Choose a continuous root $b(z)=f(z)^{1/N}$ on $\disc_{\delta}^{n-1}\times
\frac{\epsilon_0}{2}$. On each   $\pz\times\disc_{\epsilon_0}$
we find a \cont\ root $f(z)^{1/N}$ that agrees with
$b(z)$ at $z=(\pz,\frac{\epsilon_0}{2})$. Thus $f(z)^{1/N}$ is \cont\ on
$\disc_\delta^{n-1}\times\disc_{\epsilon_0}$.
Off the zero set of $f^{1/N}$, we have
$\frac{\db f^{1/N}}{f^{1/N}}=\db\log f^{1/N}=\frac{1}{N}\dff$.

Assume that $N=1$.
The local boundedness of $v$ follows from the above estimate $|g|<1/4$.
Now $$\frac{1}{v}=(1+g(z))\frac{z_n-r(\pz)}{(1+g(z))f(z)}$$
is also locally bounded,   by applying the maximal principle to
  the \holfn\ $\frac{z_n-r(\pz)}{(1+g(z))f(z)}$.
\end{proof}

\myrems{ (i)
We are not claiming that the  above $v$ (for $r_1=\cdots=r_N$)
is \cont\ in $z$, though this holds eventually.  In fact if $v$ is \cont\ and if
additionally
$\db(f|_L)=A_Lf$ on $L\cap\disc^n$ with
$\|A_L\|_{L^p(L\cap\disc^n)}\leq M$ for {\it all} complex lines $L$,  then one
can conclude that $r$ is \hol.  To see that, consider
the line   $L=\{a+tb\colon t\in\cc\}$  where $a=(a',a_n)$ is a zero of $f$ with $|a'|<\delta$ and $b$ is
a non-zero vector in $\cc^n$. We know that there is a \contfn\ $u$ on $L$ such that
$e^{-u}f^{1/N}=ve^{-u}(z_n-r(\pz))$ is holomorphic on $L$. Thus
$\lim_{t\to0}\frac{(a_n+tb_n-r(a'+tb'))v(a+tb)}{t}$
exists. Since $a_n=r(a')$ and $v(a)\neq0$,
then $\lim_{t\to0}\frac{r(a'+tb')-r(a')}{t}$ exists for all $b'$.

(ii) Let $p>2$. Assume that $f$ is continuous on $\disc^n$ and
$C^1$-smooth on $\ov\disc^n\setminus f^{-1}(0)$ and $\dd{f}{\ov z_n}=a_Lf$ on
$L\cap\disc^n$
with $\|a_L\|_{L^p(L\cap\disc^n)}\leq M$ for all $z_n$-lines $L$.
Then the $v$ (for $r_1=\cdots=r_N$)  is \cont.
To see that, we let $L=L(z)$ be the $z_n$-line passing through $z$.
 Assuming $N=1$, put
$u=A_L\ast\frac{1}{\pi z_n}$, where the convolution is on $|z_n|<1$. We already know
that $r$ is \cont.
Fix $w\in\disc^n$. We want to show that $u$ is \cont\ at $w=(\pw,w_n)$.
Split $a_L\ast\frac{1}{\pi z_n}$
into two parts:
$$
u(z)=\frac{i}{2\pi}\int_{|\zeta-r(\pw)|<\delta,|\zeta|<1}\frac{a_L(\zeta)}{z_n-\zeta}\,
d\zeta\wedge d\ov\zeta
+\frac{i}{2\pi}\int_{|\zeta-r(\pw)|>\delta,|\zeta|<1}\frac{\dd{f}{\ov
\zeta}(\pz,\zeta)}{(z_n-\zeta)f(\pz,\zeta)}
\, d\zeta\wedge d\ov\zeta.
$$
Since $\frac{1}{z}\in L^q_{loc}$ for $1\leq q<2$,  the first term tends to $0$ as  $\delta\to0$.
 For a fixed $\delta$, the second term is \cont\ in $z$ at $w$. Thus $u$
 is \cont\ on $\disc^n$. As before by  Rado's
theorem $e^{-u}f$ is \hol\ in $z_n$.
Since $e^{-u}v=\frac{e^{-u}f}{z_n-r(\pz)}$ is holomorphic in $z_n$, then
$$
(e^{-u}v)(z)=\frac{1}{2\pi i}\int_{|\zeta|=\epsilon_0}
\frac{e^{-u(\pz,\zeta)}f(\pz,\zeta)}{(\zeta-r(\pz))(\zeta-z_n)}\, d\zeta,
\quad |z_n|<\epsilon_0.
$$
Thus $e^{-u}v$ and hence $v$ is \cont.
}

\le{n1}
{Let $2<p\leq\infty$.
Let $f(z)$ be a \cont\ function on $\disc^n$, satisfying the condition (L) (with $\Omega=\disc^n$).
Assume that for each $\pz\in\disc^{n-1}$,
$f(\pz,z_n)$ has exactly one zero in $z_n=r(\pz)\in\disc$,
and its
multiplicity  is $1$. Then
$B$, the trivial extension of $\dff$ by $0$ on $f^{-1}(0)$,
 is $\db$-closed in the sense of distribution on $\disc^n$,
and $r$ is \hol.
}
\bp  $A=\dff$ is obviously $\db$-closed on the complement of
$f^{-1}(0)$. The main point in the proof consists in proving that
its trivial extension,  by 0 on $f^{-1}(0)$,  is  $\db$-closed on
$\disc^n$ (a problem that of course does not exist for $n=1$).
In case $\dff$ is bounded, this
is straightforward if one assumes (but we don't) that $f^{-1}(0)$ has a basis of
neighborhoods whose boundaries have  $(2n-1)$-Hausdorff measure
tending to 0.
The trivial extension of $A=\dff$ is precisely the $B\in L^p_{loc}(\disc^n)$ in remark 3 a).
 So $\db f=Bf$ on $\disc^n$ (not just on $\disc^n\setminus f^{-1}(0)$).

 Return to our setting.   We want to show that   $\db B=0$ near $0$.
Let $d_v(z)$ be the vertical distance from $z$ to the graph
$ z_n=r(\pz)$. So $d_v(z)=|z_n-r(\pz)|$.
Let $d(z)$ be the Euclidean distance from $z$ to the graph.
  We will prove that   $\epsilon$-neighborhoods of the graph $z_n=r(\pz)$
with respect to $d$ and $d_v$ respectively have comparable volumes. This is achieved by comparing
 $d(z)$, $d_v(z)$ with $|f(z)|$.

If   $\tilde f(z)=f(\epsilon z)$ and $\dd{f}{\ov z_n}=af$ (the one-dimensional
distribution derivative off $f^{-1}(0)$
with $z_1,\dots, z_{n-1}$ being fixed first)
we have $\dd{\tilde f}{\ov z_n}=\tilde a\tilde f$ off $\tilde
f^{-1}(0)$ with
\aln{
(\int_{|z_n|<1, \tilde f(z)\neq0}|\tilde a|^p\frac{i}{2}dz_n\wedge d\ov z_n)^{\frac{1}{p}}&
=(\int_{|z_n|<1,f(z)\neq0}|\epsilon a(\epsilon z)|^p\frac{i}{2}dz_n\wedge d\ov z_n)^{\frac{1}{p}}\\
&=\epsilon^{1-\frac{2}{p}}(\int_{|z_n|<\epsilon,f(z)\neq0}
|a(\epsilon z',z_n)|^p\frac{i}{2}dz_n\wedge d\ov z_n)^{\frac{1}{p}}.
}
Since $p>2$, by a dilation,
we may assume that the $L^p$ norms of $A_L$ on $L\cap\disc^n\setminus f^{-1}(0)$
 are small on all $z_j$-lines $L$. On
each $z_j$-line $L$, there is a \contfn\ $g$
 such that $(1+g)f$ is holomorphic on $L\cap\disc^n$ and $|g|<1/4$.
We may also assume that $|f|<1$.

Let $f(0)=0$.
Let $z\in\disc_{1/4}^n$ with  $f(z)\neq0$.

To compare $f(z)$ with $d_v(z)$, recall from \rl{2}~ii that
$$
f(z)=(z_n-r(\pz))v(z),\quad z\in\disc^n,
$$
where  $\frac{1}{v}$ is locally bounded. Restricting
  to a   smaller   polydisc if necessary, we may assume  $|v|>1/c$ on $\disc^n$ for a
fixed constant $c>0$.
  Thus $d_v(z)=|z_n-r(\pz)|\leq c|f(z)|$.

To compare $f(z)$ with $d(z)$
take $z^*\in f^{-1}(0)\cap\ov\disc_{1/2}^n$
such that
$d(z)=|z-z^*|$. Let
 $z^*=(z_1^*,\dots,z_n^*)$, $z=(z_1,\dots,z_n)$, and
 $w_j=(z_1,\dots, z_{j-1},z_j^*,\dots, z_n^*)$. Start with $w_1\in f^{-1}(0)$
 and connect $w_j$ to $w_{j+1}$ by a $z_j$-line $L_j$.
 We want to show that $|f(w_j)|\leq 8(j-1)2^jd(z)$.
 The inequality is trivial for $j=1$. Assume the inequality holds for $f(w_k)$.
  We find a function $g$, $|g|<1/4$, such that
$(1+g)f$ is holomorphic on $L_k\cap\disc^n$. By the maximum principle, we have
$$|((1+g)f)(w_{k+1})-((1+g)f)(w_k)|\leq 8|w_{k+1}-w_k|\leq 8d(z).$$
Hence
$|f(w_{k+1})|\leq 2|f(w_k)|+16 d(z)\leq 8k2^{k+1}d(z).$ Therefore $|f(z)|\leq 8n2^nd(z)$.
Thus
$$
d(z)\leq d_v(z)\leq c_1d(z), \quad z\in\disc_{1/4}^n.
$$

 Let $V_\epsilon$ be the set of points $z\in\disc_{1/4}^n$  with $d(z)<\epsilon$.
 Let $\hat V_\epsilon$ be the set of $z\in\disc^n$ with
vertical distance $d_v(z)<\epsilon$.  We have
$$
\vol(V_\epsilon)\leq\vol(\hat
V_{c_1\epsilon})=\pi(c_1\epsilon)^2\vol(\disc^{n-1})
=\tilde c\epsilon^2.
$$

Let $\chi$, $0\leq\chi\leq1$, be a smooth function such that $\chi(t)=1$ for
$t>1$, $\chi(t)=0$
for $t<1/4$, and $|\chi'|<2$. Let
$\chi_\epsilon(z)=\chi(\frac{d(z)}{\epsilon})$.
Then $|\chi_\epsilon(z')-\chi_\epsilon(z)|\leq \frac{2|z'-z|}{\epsilon}$. Let
$\chi_\epsilon^*$ be a regularization of $\chi_\epsilon$
  such that $\chi_\epsilon^*$
vanishes on $V_{\epsilon/8}\cap\disc_{1/2}^n$ and $\chi_\epsilon^*=1$ on
$\disc_{1/2}^n
\setminus V_{3\epsilon/2}$.
Note that $0\leq\chi_\epsilon\leq1$ and
$|\nabla\chi_\epsilon^*|\leq2/\epsilon$.

Let $\varphi$ be a smooth $(n,n-2)$-form
supported in $\disc_{1/4}^n$.
Let $\frac{1}{p}+\frac{1}{q}=1$.  Since $p\geq 2$ and $1\leq q\leq 2$,
  we have
\aln{
&\Bigl|\int_{\disc^n}\chi_\epsilon^*A\wedge\db\varphi\Bigr|
  =\Bigl|\int_{\disc^n}A\wedge\db\chi_\epsilon^*\wedge\varphi\Bigr|\\
&\qquad \leq c_0
\|\varphi\|_{L^\infty}\Bigl(\int_{V_{3\epsilon/2}}\left|A\right|^p\,
dV\Bigr)^{1/p}
\Bigl(\int_{V_{3\epsilon/2}}|\db\chi_\epsilon^*|^q\,dV \Bigr)^{1/q}\\
&\qquad \leq c_1
\|\varphi\|_{L^\infty}\Bigl(\int_{V_{3\epsilon/2}}\left|A\right|^p\,
dV\Bigr)^{1/p}
\epsilon^{\frac{2}{q}-1}
\to 0,\quad\text{as $\epsilon\to0$}.}

\

Therefore,
$B$, the trivial extension of $\dff$ by $0$ on $f^{-1}(0)$,  is $\db$-closed on
$\disc_{1/4}^n$ (and hence on $\disc^n$).

\

We are ready to show that $r$ is holomorphic.
By Proposition C and remarks~3~a)-b),
the uniform $L^p$ bound for $A_L$
  in $\db (f|_{L\cap\disc^n})=A_Lf$, implies that $\dd{f}{\ov
z_j}=B_jf$ on $\disc^n$.
Moreover, for all $z_j$-lines
$L$ and the $\epsilon$-neighborhood $L_\epsilon$ of $L$ in $\cc^n$ with  the sup-norm,
  $\|B_j|_{L_\epsilon\cap D}\|_{L^p}\leq (\pi\epsilon^2)^{\frac{n-1}{p}}M$.
By remark~3~d), all solutions $u$ to  $\db u=B$, which exist, are in $C_{loc}^\alpha(\disc^n)$.
Thus $e^{-u}f$ is \hol\ on
$\disc^n\setminus
f^{-1}(0)$, and hence on $\disc^n$ by Rado's theorem. Since the zero set of
$e^{-u}f$
is given by the graph  $z_n=r(\pz)$, then $r$ is \hol.
\ep


\myrems{
(a) The above estimates on two distances $d(z)$, $d_v(z)$ are crucial. One can
obtain $|f(z)|\leq c\|A_{L_j}\|_{L^p}d(z)^{1-\frac{2}{p}}$ from equations
$\db (f|_{L_j})=A_{L_j}f$ by the
well-known estimates stated at the beginning of Section 2, without using an integrating factor.
When $p\geq 3+\sqrt 5$ (i.e $2(1-\frac{1}{p})(1-\frac{2}{p})\geq1$)
that estimate  insures   the above lemma,   since as before
\aln{
&\Bigl|\int_{\disc^n}\chi_\epsilon^*A\wedge\db\varphi\Bigr|\leq c_0
\|\varphi\|_{L^\infty}\Bigl(\int_{V_{3\epsilon/2}}\left|A\right|^p\,
dV\Bigr)^{1/p}
\Bigl(\int_{V_{3\epsilon/2}}|\db\chi_\epsilon^*|^q\,dV \Bigr)^{1/q}\\
&\qquad \leq c_0
\|\varphi\|_{L^\infty}\Bigl(\int_{V_{3\epsilon/2}}\left|A\right|^p\,
dV\Bigr)^{1/p}
\Bigl(\int_{{\hat V}_{c_1\epsilon^{1-\frac{2}{p}}}}|\db\chi_\epsilon^*|^q\,dV \Bigr)^{1/q}\\
&\qquad \leq c_2
\|\varphi\|_{L^\infty}\Bigl(\int_{V_{3\epsilon/2}}\left|A\right|^p\,
dV\Bigr)^{1/p}
\epsilon^{2(1-\frac{1}{p})(1-\frac{2}{p})-1}
\to 0,\quad\text{as $\epsilon\to0$}.}
(b) To obtain the $\db$-closedness of the trivial extension
of $\dff$ on $\disc^n$,  it is crucial to assume that
$\db(f|_L)=A_Lf$ on $L\cap \disc^n$
with $\|A_L\|_{L^p}\leq M$ for all $z_j$-lines and for all $j=1,\dots, n$.
Take $n=2$ as an example. We do not know if the trivial extension is $\db$-closed if
 the $\db$-inequality is assumed on
 all $z_1$-lines and  if the $\db$-inequality on all $z_2$-lines is replaced
  by $\dff\in L^p(\disc^2)$ (for any $2<p<\infty$).
The latter insures  $\db(f|_L)=A_Lf$ on $L\cap\disc^2$ with $A_L\in L^p$
for almost all $z_2$-lines $L$, but the $L^p$ norms of
$A_L$ might not be bounded as $L$ varies.

\sk
{\bf 3.3. End of proof of Theorem A.}

\sk

  The theorem is local. We may assume
that $f$ is defined on $\disc^n$. For each $z\in\disc^n$,
 denote by $L_j(z)$ the $z_j$-line that is parallel to the $z_j$-axis.
 By one-dimensional
result, we know that   for each $p\in\disc^n$ and each complex
line $L=L_j(p)$, there is a \contfn\ $u$ so that
$e^{-u}f|_L$ is holomorphic on $L\cap\Omega$. Define the order of
vanishing of $f|_L$   at $p$ to the vanishing order of
$e^{-u}f|_L$ at $p$, which is independent of the choice of $u$.
  Let $N_p=\infty$ if $f|_{L_j(p)}\equiv0$ for all $j$;
otherwise   let $N_p$ be the smallest integer such that
for   some   $L=L_j(p)$, $f|_L$ vanishes to order
$N_p$ at $p$.

\bigskip
Assume that $f(0)=0$. Since we consider only lines parallel to the coordinate axes
it is possible that $N_0$ is infinite although $f$ is not identically 0.
By permuting the coordinates, we may assume that
  $z_n=0$ is a zero of $f(0,z_n)$ of order
$N_0$.

\medskip
(i) Case  $N_0< +\infty$.

\medskip When $N_0=1$, near the origin,
$f^{-1}(0)$ is a smooth complex
hypersurface by \rl{2}~i and \rl{n1}.

Assume that $V$ is a
complex variety near $p$ if $N_p<N$.

\medskip
If  $r_1=\cdots=r_N$ then by \rl{2} ii, there is a \cont\ $N$-th root
$f^{1/N}(z)=(z_n-r(\pz))  v(z)$, which still satisfies the $\db$-differential
inequality. We have
$N_0(f^{1/N})=N_0(f)/N=1$. Hence $f^{-1}(0)$ is
a smooth complex hypersurface near the origin.

\medskip
Suppose $N_0=N\geq2$. Without loss of generality, we may assume by
\rl{2} i  that the zero set of $f$   in $\disc^n$
  is given by $z_n=r_j(\pz),\ j=1,\dots, N$, counting multiplicity.

We already proved the assertion when  $r_1=\cdots=r_N$.
So we may assume that  not all of
  $r_1(a),\dots, r_N(a)$ are the same for some $a$. Let $k$ be the largest
integer such that
  $r_1(\pz),\dots, r_N(\pz)$ have $k$ distinct elements for some $\pz$.
  Define $$h(\pz)=\Pi_{1\leq \alpha\neq \beta\leq k}
  (r_{j_\alpha}(\pz)-r_{j_\beta}(\pz)),$$ when
 $r_{j_1}(\pz),\dots, r_{j_k}(\pz)$  are distinct.
 Rename the distinct $k$ elements
 by $r_1^*(\pz), \dots,$ $ r_k^*(\pz)$.
  For clarity, we do not define $r_1^*(\pz),
 \dots, r_k^*(\pz)$ when $\{r_1(\pz),\dots, r_N(\pz)\}$ has less than $k$
  distinct elements, in which case we set $h(\pz)=0$.

We first want to show that $h$ is holomorphic away from $h^{-1}(0)$. Assume that $h(a)\neq0$.
Since $k\geq2$ then $N_{(a,r_{j}^*(a))}<N$.
Thus $f^{-1}(0)$ is a complex variety near each $(a,r_{j}^*(a))$.
By \rl{2} i, for $\pz$ close to   $a$, $f(\pz,z_n)$ has at least one zero near each
 $r_{j}^*(a)$.
 Therefore, by the definition of $k$, $f(\pz,z_n)$ has
exactly one  zero near each $r_{j}^*(a)$
 for $\pz$ sufficiently close to $a$.   Thus near each $(a,r_{j}^*(a))$
 the complex variety $f^{-1}(0)$ must be smooth, and near $a$  we redefine  $r_{j}^*(z')$ such that they
 become holomorphic
 in $z'$.
In particular $h$ is holomorphic away from its zero set.

  Next we want to show that $h$ is \cont.
  Fix $a\in\disc^{n-1}$ with $h(a)=0$.  Given a sequence $a^u$ approaching to $a$ as $u\to\infty$
   we want to show
  that $\lim_{u\to\infty}h(a^u)=0$. Without loss of generality we may assume that $h(a^u)\neq0$
  for all $u$.
     Let $b_1,\dots,
b_d$ be distinct zeros of $f(a,z_n)$ in $\disc\ni z_n$.
  Let $m_j$ be the multiplicity of the zero $z_n=b_j$ of $f(a,z_n)$. Put
 $s=\frac{1}{4}\min\{|b_j-b_k|, 1-|b_k|\colon 1\leq k\neq j\leq d\}$.  By \rl{2} i,
we can choose $t>0$ so
that
the zero set of $f$ in $(a+\disc_t^{n-1})\times (b_j+\disc_s)$ is
given by $z_n=w_{j\alpha}(\pz)$, $1\leq\alpha\leq m_j$; moreover
$\lim_{\pz\to a}w_{j\alpha}(\pz)=b_j$. Therefore the limit of any
convergent subsequence  $(r_{1}^*(a^{u_j}),\dots, r_{k}^*(a^{u_j})), j=1,2,\dots$
must have the form $(b_{l_1},\dots,b_{l_k})$. Since $\{b_{l_1},\dots,b_{l_k}\}\subset\{b_1,\dots,
b_d\}$  and $d<k$, we conclude that $\lim_{j\to\infty}h(a^{u_j})
=\prod_{1\leq \alpha\neq\beta\leq k}(b_{l_\alpha}-b_{l_\beta})=0$.
Therefore $h$ is \cont\ and hence \hol\ by Rado's theorem.

 Away from the zero set of $h$, the symmetric polynomials of   $r_{1}^*(\pz),\dots, r_{k}^*(\pz)$
are holomorphic in $\pz$ (in fact we proved that away from the zero set,
$r_{1}^*(\pz), \dots, r_{k}^*(\pz)$ can be locally rearranged
to become holomorphic). Since the symmetric polynomials are bounded,
  they extend holomorphically to   $\disc^{n-1}$, by the removable
singularity
  theorem. By the removable
singularity
  theorem again,
   $P(z)=(z_n-r_{1}^*(\pz))\cdots(z_n-r_{k}^*(\pz))$, defined for $h(\pz)\neq0$ only,
    extends   \hol ally to $\disc^n$. Now we want to find a neighborhood $U$
    of  the origin such that  $P^{-1}(0)\cap U=
    f^{-1}(0)\cap U$.
 By definition  $P^{-1}(0)\cap\disc^n\setminus h^{-1}(0)=
   f^{-1}(0)\cap\disc^n\setminus h^{-1}(0)$. Since $h\not\equiv0$
   the closure of $f^{-1}(0)\setminus h^{-1}(0)$ in $\disc^{n}$ is $f^{-1}(0)$ by \rl{2} i.
    Since $P(0,z_n)\not\equiv0$ there is a polydisc $U=\disc_\delta^{n-1}\times\disc_\epsilon$
    such that  $P^{-1}(0)\cap U$ is given by
    a branched-covering over $\disc_\delta^{n-1}$. So $\ov {(P^{-1}(0)\setminus h^{-1}(0))}\cap U
    =P^{-1}(0)\cap U$.
Therefore  $f^{-1}(0)\cap U=P^{-1}(0)\cap U$, which is a complex variety.

\bigskip
(ii) Integrating factor and case $N_0=+\infty$.

\medskip
We now go back to the main question which is the $\overline\partial$-closedness of $B$
(, which was already central in \rl{n1}).

\medskip
We shall use the following result of Demailly ([D], Lemma 6.9), whose proof is similar to the
argument in \rl{n1}, with an induction on the dimension of the analytic set.
Let $\Omega$ be an open set in $\cc^n$, and let $E$ be an analytic subset of
$\Omega$ of dimension $<n$. Let $g$ be a   $(0,1)$-form   defined
on $\Omega$ with $L^2$ coefficients. If $g$ is $\overline\partial$-closed
on $\Omega \setminus E$ then $g$ is $\overline\partial$-closed on $\Omega$.

\medskip
We now assume that Theorem A is proved in dimension $n-1$ and we want to establish it
in dimension $n$. Assume that there exists $a=('a,a_n)\in \disc^n$ with
$f(a)\neq 0$. Let
$X=\{ 'z\in D^{n-1}:f('z,a_n)=0\}$. By the induction hypothesis, this is a proper analytic subset of
$\disc^{n-1}$.
If $p=('p,p_n)$ is such that $N(p)=+\infty$ then $f('p,.)\equiv 0$, so $'p\in X$.

\smallskip
Set $Z=X\times \disc$. This is a proper analytic subset of $\disc^n$, and
if   $p\not\in Z$  then
$N(p)<+\infty$.

\bigskip
Let $B$ be the trivial extension of ${\overline\partial f \over f}$, as considered earlier.
  We already know that $B\in L_{(0,1)}^p(\disc^n)$.
We claim that $B$ is $\overline\partial$-closed on $\disc^n$. Applying Demailly's result
to $\Omega = \disc^n\setminus Z$, it follows from the analyticity of $f^{-1}(0)$
(at points where $N$ is finite)
shown in
(i) that $B$ is $\overline\partial$-closed on $\disc^n\setminus Z$. Applying Demailly's result
again, it follows that $B$ is $\overline\partial$-closed on $\disc^n$, as desired.
As at the end of the proof of \rl{n1},   for some locally defined $u\in C_{loc}^\alpha$,
$e^{-u}f$ is  holomorphic.

The proof of Theorem A is complete.

\sk
{\bf 3.4. Vector-valued $f$ defined on $\cc^n$ when  $f^{-1}(0)$ is real analytic.}

\sk

We turn to the case that $f$ is vector-valued and defined on a domain $\Omega\subset\cc^n$
with $n>1$. Here our result is far from complete. We will treat only the case that $f^{-1}(0)$
is already real analytic. To show the complex analyticity, we will impose condition on $\db(f_\gamma)$
for all germs of complex curve, not just  lines. But we will not need the uniform bound on $L^p$-norms.

\

First we recall some basic results on real analytic sets, which can be found in~[N].

Let  $V_0$ be an irreducible germ  of  real analytic set  at $0\in\rr^m$ of dimension $k$. (We will take
$\rr^m=\rr^{2n}\subset\cc^{2n}$ soon.)
Then
there is an open subset $D$
of $\rr^m$ and a closed real analytic set $V$ in $D$ which represents the germ, with $\dim V_x\leq k$
for all $x\in V$. Denote by $V^*$ the set of points in $V$ at which $V$ is a smooth submanifold of dimension
$k$.
We also need the complexification $\tilde V_0$ of the germ $V_0$:
 there exists a unique irreducible germ $\tilde V_0$ of complex variety
in $\cc^m$  with $\tilde V_0\cap\rr^m=V_0$
 ([N], Proposition 1, p.~91). In particular, one can choose an open subset $\tilde D$
of $\cc^m$ such that $\tilde V_0$ is represented by an irreducible closed complex variety $\tilde V\subset\tilde D$
of pure dimension $k$.

We need the following.
\le{df}{
Let $V_0$ be an irreducible germ of real analytic variety at $0\in\cc^n$ of dimension $k$,
represented by  a closed real analytic set $V$ in $U\ni0$ of the same dimension.
If $V^*$ is a complex submanifold of $\cc^n$, then
$V$ is a complex variety at $0$.
}
\bp
It suffices to find a germ of complex variety $\hat V$ at $0$ that is contained in
$V$ and has the same real dimension as $V$.
Then the irreducibility of $V$ implies that  two germs  must
agree
([N], Proposition 7, p.~41).
We will adapt an argument in [D-F] to construct $\hat V$. (The argument in [D-F]
is for
a smooth real hypersurface $V$ in $\cc^n$.)

Choose
 a polydisc $D\subset U$, centered at $0$, such that $V\cap D$ is the zero set of a real
function $r(z,\ov z)$ which is a convergent power series on $D\times D$.
 Thus $Q_w\df\{z\in D\colon r(z,\ov w)=0\}$ is a complex variety in $D$ for each
$w\in D$. Choose a domain $\tilde D_1\subset\cc^{2n}$
and an irreducible  complex variety $\tilde V$ in $\tilde D_1$ such that $\tilde V_0$
is the complexification of $V_0$.  We may assume that $\tilde D_1\cap\rr^{2n}\subset D$ and $\tilde V\cap\rtn
=V\cap\tilde D_1$.

We want to show that
 $Q_w\supset V^*\cap \tilde D_1$ for $w\in V^*\cap \tilde D_1$. To that end,
 choose a bi\holm\  $\varphi\colon
 \disc^{k}\to W\subset V^*$ with $\varphi(0)=w$. Then $r(\varphi(t),\ov{\varphi(t)})=0$ for $t\in\disc^n$.
 Thus $r(\varphi(t),\ov{\varphi(0)})=0$, i.e $r(\cdot,w)$
 vanishes on some open subset $W$ of $V^*$.
 Note that $W\subset\cc^n=\rr^{2n}+i0$ is
 embedded in $\cc^{2n}$. Thus $W$ is a totally real submanifold in $\tilde V$ of dimension $k$.
 As a function in $(x,y)\in\tilde V^*$,
 $r(x+iy,w)$ vanishes on an open subset of $\tilde V^*$. Since $\tilde V^*$ is connected,
 then $r(x+iy,w)$ vanishes on $\tilde V$, which implies that $r(z,w)$ vanishes for $z\in
  V^*\cap\tilde D_1(\subset\tilde V\cap\rr^{2n})$.

 We just proved that $S\df\cap_{w\in V^*\cap \tilde D_1}Q_w\supset V^*\cap \tilde D_1$.
 For $w\in V^*\cap \tilde D_1$
  we have $S\subset Q_w$, i.e $z\in Q_w\ (\Leftrightarrow w\in Q_z)$ for $z\in S$.
  Thus $V^*\cap \tilde D_1\subset Q_z$ for  $z\in S$.
 Now
 $$
 V^*\cap\tilde D_1\subset \cap_{w\in S}Q_w\df\hat V.
$$
Fix $z\in\hat V\subset S$. We have $z\in S$ and get $z\in (\hat V=)\cap_{w\in S}Q_w\subset Q_z$.
Now $z\in Q_z$ implies that $r(z,\ov z)=0$, i.e $z\in V$. Finally, $V^*\cap\tilde D_1\subset\hat V\subset V$.
Since $V^*$ has pure complex dimension $k$ and $0\in\ov V^*$ then $\hat V$ has real dimension at least $2k$,
 which is already the dimension of $V$. Since $V$ is irreducible then $V=\hat V$ as germs at $0$.
\ep

As a consequence of \rl{1} and \rl{df}, we have
\pr{an}{
Let $2<p\leq\infty$.
Let $f=(f_1,f_2,\dots, f_m)$,
where $f_j$ are \contfn s on a  domain
$\Omega\subset\cc^n$, and let $\Omega^*=\Omega\setminus f^{-1}(0)$.
Assume
that   for   each smooth \holc\ $\gamma$ in $\Omega$ (not necessarily closed in $\Omega$),
$\db (f|_\gamma)=f A_\gamma $
holds
in the distribution sense
on $\gamma\cap\Omega^*$,
where $A_\gamma$ is an $m\times m$ matrix of $(0,1)$-forms whose coefficients are
 in $L^p(\gamma)$. If
$f^{-1}(0)$
is a   real analytic variety in $\Omega$, then it is  a complex
variety.
}
\bp Assume that $f(0)=0$. We want to show that $f^{-1}(0)$ is a complex variety
near the origin.

Decompose the germ of
$f^{-1}(0)$ at $0$ as $V_0^1\cup\dots\cup V_0^k$, where $V_0^j$ are germs of
irreducible real analytic sets at $0$ and  $V_0^j$ is not contained in $\cup_{k\neq j}V_0^k$.
(For real analytic sets, the irreducible decomposition may
exist at the germ level only!) Choose an open neighborhood $D$ of $0$ and  a
closed real analytic set $V_j$ in $D$ such that $V_0^j$ is the germ of $V_j$ at $0$.
  Since $f^{-1}(0)$ and $\cup V_j$ represent the same germ at $0$ there exists
an open set $D\ni 0$ such that
 $V\df f^{-1}(0)\cap D=\cup V_j\cap D$.
Let $p_j$ be the dimension of $V_0^j$. We may assume that $V_j$ has dimension $p_j$ too.
Since $\dim\cup_{k\neq j}V_0^k\cap V_0^j<\dim V_0^j$, we may choose $D$ so small that
 $\dim\cup_{k\neq j}V_k\cap V_j<\dim V_j$. Let $V_j^*$ be the set of points in $V_j$
 near which $V_j$ is a real submanifold of dimension $p_j$. Put $\tilde V_j=
 V_j^*\setminus\cup_{k\neq j}V_k$. Each tangent vector in $T_x\tilde V_j$
 is the tangent vector of some real analytic curve $\gamma$
 in $\tilde V_j$. Let $\tilde \gamma$ be the complexification of $\gamma$. By assumption,
 $f^{-1}(0)\cap\tilde\gamma$ is a complex variety in $\tilde \gamma$. Since $f^{-1}(0)\cap\tilde\gamma$
 is not isolated, then $\tilde \gamma$ is contained in $f^{-1}(0)$ (one may assume that $\tilde\gamma$
 is connected). Therefore, $T_x\tilde V_j$ is complex linear subspace
    of $T_x\cc^n$. Since $\tilde V_j$ is dense in $V_j^*$, then
     $V_j^*$ is a complex submanifold of $\cc^n$. By the previous lemma, we know
  that each $V_j$ is a complex variety.
\ep

\setcounter{thm}{0}\setcounter{equation}{0}
\section{Remarks on some uniqueness or finite type
results for \jholc s}

In the theory of \jholc s inequalities of the type
$$
|\db f|\leq c(z)|f|  \eqno{(T 1)}
$$
occur, as well as inequality of the type
$$
|\db f|\leq c(f)|\partial f|.  \eqno{(T 2)}
$$
The more general inequality $|\db f|\leq c(z,f)|\partial f|$
can,  of course,   be
reduced to (T2) by considering $(z, f(z))$
instead of $f(z)$.

If one has an $L^p$-bound for $\partial f$ and an estimate $|c(f)|\leq
|f|$, reversing roles, one can reduce (T2) to (T1).

An example is the following: Assume that $v_1\colon \disc\to(\rtn, J)$
is an embedded \jholc. After changing variables we assume that $v_1(z)
=(z,0,\cdots, 0)\in\cc^n\simeq\rtn$, and that $J(z,0,\cdots,0)=J_{st}$
(the standard complex structure on $\cc^n$). Let
$v_2\colon \disc\to(\rtn, J)$ be \cont.   Let $E$ be a closed subset of $\disc$
with $0\in E$.
Assume that on $E$ $v_2(z)=v_1(z)=(z,0,\cdots,0)$ and that $v_2$ is \jhol\
on $\disc\setminus E$ (but make no a-priori  assumption of J-holomorphicity of $v_2$
at the boundary points of $E$). It follows from  \rl{1} in Part One  that  if $\nabla v_2
\in L^p$ for some $p>2$ then near $0$, $v_2\equiv v_1$ or $v_2-v_1$ vanishes to
finite order only. Indeed, the equation for J-holomorphy is $\dd{v}{\ov z}=Q(v)
\ov{\dd{v}{z}}$, and $Q(z,0,\dots,0)\equiv0$, hence $Q(v_1)\equiv0$ (with
$Q\in M_{n,n}(\cc)$ and of the smoothness of $J$).
So
\aln{
|\dd{(v_2-v_1)}{\ov z}|&=|\dd{v_2}{\ov z}|=|Q(v_2)\ov{\dd{v_2}{z}}|
=|\dd{v_2}{z}| |Q(v_2)-Q(v_1)|\\ &\leq c(z)|v_2-v_1|,\qquad \text{with $c\in
L^p$.}}

However   the above result (that $v_2\equiv v_1$ or $v_2-v_1$ vanishes to finite order only)
follows immediately from the Corollary to
Theorem C of Part Two (applied to $(z,v_j)$ rather than $v_j(z)$).
Indeed this Corollary implies that in fact $v_2$ is $J$-holomorphic
(i.e. the exceptional set $E$ can be removed).
Then,
the result is well-known.

\

\begin{center}
\bf  Part Two:
Removable singularities for \jholm s and Rado's Theorem
\end{center}

\sk
{\bf 1. Introduction}

\

Let $\Omega$ be an open set in $\cc$. A subset $E\subset\Omega$ is said to be
a {\it
  polar} set if for any $z\in\Omega$ there exists a subharmonic function
$\lambda$ defined on a neighborhood $V$ of $z$ (not identical to $-\infty$
near $z$) such that $E\cap V\subset\lambda^{-1}(-\infty)$. If $E$ is closed
in $V$ one can take $\lambda$ to be \cont\ as a map from $V$ onto
$\{-\infty\}\cup\rr$
(and $E\cap V=\lambda^{-1}(-\infty)$).

If $E$ is a polar subset of $\Omega$, then there exists a \subfn\ $\mu$ on
$\Omega$
$\mu\not\equiv-\infty$ such that $E\subset\mu^{-1}(-\infty)$.
 Indeed local polarity implies global polarity, and this is an easy result
unlike Josefson's theorem on pluripolarity.
  (Hint:
Take a locally finite covering of $\Omega$ and $\lambda_j$ corresponding
\subfn s and solve $\Delta\lambda=\sum_j\Delta\lambda_j$.)

In  Part Two
 all \acs s are of class $C^1$, and all
  \jholc s $u$ are of class $C^1$ at least. So $u$ are of class $C^{1,\alpha}$
  for   any   $\alpha<1$ (See [I-R]).

\sk
{\bf Theorem B}. {\em Let $\Omega$ be an open subset in $\cc$ and let $E$ be a
closed polar
subset
of $\Omega$. Let $u$ be a \contm\ from $\Omega$ into an \acm\ $(M,J)$ with $J$
of class $C^1$. If $u$ is \jhol\ on $\Omega\setminus E$ then it is \jhol\ on
$\Omega$.}

The next theorem is
also about removable singularities but this time in terms of
the target space. It is a theorem in the style of the classical theorem
of Rado that states that a {\em \cont} function that is \hol\ off its zero
set is \hol.

\sk
{\bf Definitions}. 1) A closed subset $\mathcal C$ of an \acm\ $(M,J)$ will be
said
to be a {\it Rado subset}
of $(M,J)$ if and only if the following holds:
If $u$ is a \contm\ from a connected open subset $\Omega$ of $\cc$
into $M$ such that $u$ is \jhol\ at any point $z$ such that
$u(z)\not\in\mathcal C$
then $u$ is \jhol\ on  $\Omega$   or $u(\Omega)\subset\sC$.

2) A subset $L$ of $(M,J)$ is said to be {\em locally exact} $J$-pluripolar, if
and only
if for any $q\in L$ there exits a neighborhood $V$ of $q$ and a \jholfn\ $\rho$
defined on $V$, \cont\ as a map from $V$ into $\{-\infty\}\cup\rr$, such that
  $\rho\not\equiv-\infty$ near $q$ and
$L\cap V=\rho^{-1}(-\infty)$. (The notion of {\em $J$-plurisubharmonicity} will
be
recalled later.)

\sk
{\bf Theorem C}. {\em Let $\sC$ be a closed subset of an \acm\ $(M,J)$.
Assume that there is a discrete subset   $S$ of $\sC$   such that $\sC\setminus
  S$
is locally exact $J$-pluripolar then $\sC$ is a Rado subset.}

\sk
{\bf Corollary}. {\em Let $(M,J)$ be an \acm\ of class $C^1$.
The proper image of an open subset of $\cc$ under a  \jholm\
of class $C^2$ is a Rado subset. Discrete subsets of $M$ are Rado subsets.}

\

By the proper image of an open subset $\Omega$ of $\cc$ under a \jholm\ $u$, we mean that
the map $u\colon \Omega\to    M$ is proper.

\

  Proposition~B, which
  asserts the existence of a small neighborhood $\omega\subset\disc$ of the origin
such that either the $v$ is \jhol\ on $\omega$
or $v(\omega)$ is contained in $ u(\disc)$,
 follows from the above Corollary.  In general,
one cannot take $\omega$ to be $\disc$.

\

Two results on $J$-pluripolarity will be used. The first one will allow us
to prove that Theorem C follows from Theorem B. The second one shows that the
corollary follows immediately from Theorem C.

a) If $p\in(M,J)$ there exists a  $J$-pluri\subfn\ $\rho$
defined on a neighborhood of $p$ and \cont\ away from $p$
such that $\rho^{-1}(-\infty)=\{p\}$. This
is a local question. We can assume that $M=\rtn\simeq\cc^n, p=0$, and
$J(0)=J_{st}$.
Then for $A>0$ sufficiently large one can take $\rho(Z)=Log|Z|+A|Z|$ (Chirka,
see a proof in [I-R] Lemma 1.4 page 2401).

b) If $\Sigma$ is a (germ of) embedded \jhol\ $C^2$ disc and $q\in\Sigma$ there
exists
a $J$-pluri\subfn\ $\rho$ defined  in a neighborhood $V$ of $q$,
which  is \contm\ from $V$ into $\{-\infty\}\cup \rr$,  such that
$\Sigma\cap V=\rho^{-1}(-\infty)$ ([R]).

\bigskip

Recall that a function $\lambda$ defined on an open set of $(M,J)$ is said
to be \jpsh\ if and only if $\lambda$ is upper semi\cont\ and its
restriction to any \jholc\ is \sub\ (i.e if $u\colon \disc\to(M,J)$ is
\jhol\ then $\lambda\circ u$ is \sub). For a smooth function $\lambda$
the condition is that for any tangent vector $T$ of $M$ at a point $q\in M$:
 $(dd_J^c\lambda )_q(T,J(q)T)\geq0$.
See Corollary 1.1 in [I-R].

We shall repeatedly use the crucial formula: If
$u$
is a $C^1$ \jholm\ and $\lambda$ is a $C^2$ function then
$$
\Delta(\lambda\circ u)(z)=[d\,d_J^c\lambda]_{u(z)}\left(\dd{u}{x}(z),
J(u(z))\dd{u}{x}(z)\right)
\eqno{(*)}
$$
(formula (1.2) in [I-R] page 2399), where $d_J^c\lambda(Y)=-d\lambda(JY)$.
In [I-R] there is no claim of originality for the above   formula,   and the reference is given
for the
convenience of the reader.

\sk
{\bf Remarks.} 1) We do not know whether \jholc s with cusps are (locally)
$J$-pluripolar.
It would   simplify  our proof of Theorem C.

2) Theorem B is a generalization of the proof of removal  of isolated
singularities
for \jholm s that extend \cont ly. Its proof borrows from the known proof
of the removal of isolated singularities, although some of arguments are
slightly different from the usual ones even in that case.

\bigskip
\noindent
{\bf 2. Proof of Theorem B}

\sk {\bf 2.1 Preliminaries.}

\sk

There is no claim of originality for the following lemmas. The first one
is just the case for the classical Rado theorem,   and   the second one is a version
of
Rado's theorem. The third one is certainly very classical.

\sk
{\bf Lemma 1.} {\em Let $\Omega$ be an open set in $\cc$ and let
$E$ be a closed polar subset of $\Omega$. If $v$ is a \contfn\ on $\Omega$
that is \sub\ on $\Omega\setminus E$,  then $v$ is \sub\ on $\Omega$.}

\noindent
{\bf Proof.} Let $\rho$ be a \subfn\ on $\Omega$ that is $-\infty$
 exactly  on $E$ and not identical $-\infty$ on any open set. Then
$v=[\lim_{\epsilon\to0^+}
(v+\epsilon\rho)]^*$, where $*$ denotes upper semi\cont\ regularization.


\sk
{\bf Lemma 2.} {\em Let $\Omega$ and $E$ be as in Lemma 1.
If $u$ is a bounded harmonic function defined on $\Omega\setminus E$
then $u$ extends to a harmonic function on   $\Omega$.}

\noindent
{\bf Proof.} Set $u^*(z)=\limsup_{\zeta\to z} u(\zeta)$ for $\zeta\in \Omega$.
Then $u^*$ is \sub\ since $u^*=[\lim_{\epsilon\to0^+}
(u+\epsilon\rho)]^*$. Also $(-u)^*$ defined similarly is \sub. At any
$z\in\Omega
\setminus E$, $u^*(z)=u(z)$ and $(-u)^*(z)=-u(z)$. So by the mean value
property for $u^*$ and $(-u)^*$,
$u$ is obtained on any disc relatively compact
in $\Omega$ by integration of $u$ against the Poisson kernel on its boundary
(whose intersection with $\rho^{-1}(-\infty)$ has zero measure).

\sk
{\bf Lemma 3.} {\em Let $\Omega$ be a domain in $\cc$ and let $E$ be
a closed polar subset of $\Omega$. Let $g$ be a bounded \subfn\ on
$\Omega$.   If $\mu = \Delta g$ (it is a positive measure), then
 $\mu(E)=0$.
}

\noindent
{\bf Proof.} The question is local, so we may assume that $\Omega=\disc$
and $\mu$ is a finite measure with compact support in $(\ov\disc)$.
Then $g=\mu*N+h$, where $N=\frac{Log|z|}{2\pi}$   is the
Newtonian potential and $h$ is harmonic on $\disc$. Let $\mu=\mu_E+\mu_{E^c}$
where $\mu_E$ is the restriction of the measure $\mu$ on $E$. We have to show
that
$\mu_E=0$. Since $g$ is bounded $\mu*N$ must be locally bounded on $\disc$.
Since both $\mu_E$ and $\mu_{E^c}$ are positive (and only large positive values
can occur to $-Log|z|$) both $\mu_E*N$ and $\mu_{E^c}*N$ (a priori in
$L^1(\disc)$)
must be locally bounded. But $\mu_E*N$ is harmonic off $E$. So by Lemma 2,
$\mu_E
*N$ extends to a harmonic function and so $\mu_E=\Delta(\mu_E*N)=0$.

\bigskip
\noindent
{\bf 2.2. The proof of Theorem B.}

\medskip

The question is local. We can assume that $M=\rtn$,  $0\in\Omega$,
$u(0)=0\in\rtn$
and in $\rtn$ $J(0)=J_{st}$. We want to prove $J$-\hol ity of $u$
at $0$.

\sk
{\bf Step 1.} $\nabla u\in L^1_{loc}(\Omega)$. Here $\nabla u$ means the
distributional gradient of $u$ on $\Omega$, not only its restriction to
$\Omega\setminus E$. We identify $\rtn$ with coordinates $(x_1,y_1,
\dots,x_n,y_n)$ with $\cc^n$ via $z_j=x_j+iy_j$, and put $|Z|^2=\sum_{j=1}^n|z_j|^2$.
Note that $|Z|^2$ is strictly \jpsh\ near $0$ (see [I-R], Lemma~1.3 page 2400).
Let $\varphi$ be either the function $\RE z_j$ or $\IM z_j$ for some $j\in\{1,
\dots, n\}$. Then $|Z|^2$ is   \jpsh\ near $0$ and so is $\varphi+K|Z|^2$
if $K>0$ is sufficiently large.

In a neighborhood of $0\in\cc$, $|Z|^2\circ u$ and $(\varphi+K|Z|^2)\circ u$
are \sub\ by Lemma~1 since off $E$ $u$ is \jhol. Therefore $\Delta((\varphi+
K|Z|^2)\circ u)$ and $\Delta(|Z|^2\circ u)$ both are positive measures,
%
and $d\mu=\Delta(\varphi\circ u)$ is  a  (locally finite) measure. Hence
$\nabla(\varphi\circ u)\in L_{loc}^1(\Omega)$,
since $\varphi\circ u$ is equal to $\Delta(\varphi\circ
u)|_{\disc}*\frac{Log|z|}{2\pi}$
modulo smooth functions and the distributional gradient
$\nabla(d\mu*\frac{Log|z|}{2\pi})=(\nabla \frac{Log|z|}{2\pi})\ast d\mu\in L_{loc}^1(\Omega)$.

The last assertion is true for each coordinate function $\varphi$, so $\nabla u
\in L_{loc}^1(\Omega)$ as claimed.

\sk
{\bf Step 2.} Now that we know that $\nabla(\varphi\circ u)\in L_{loc}^1$ we
can
show that $\nabla u\in L^2_{loc}(\Omega)$.

Off $E$ (which has Lebesgue measure $0$ since it is polar),
$$
\Delta(\lambda\circ
u)(z)=[d\,d_J^c\lambda ]_{u(z)}\left(\dd{u}{x}(z),J(u(z))\dd{u}{x}(z)\right).
$$
For $p$ in a neighborhood of $0$ we have for any tangent vector $T$ at $p$
$$
d\, d_J^c|Z|^2(T,J(p)T)\geq C|T|^2.
$$
So $|\dd{u}{x}|^2\leq\frac{1}{C}\Delta(|Z|^2\circ u)$.
Since $\dd{u}{y}=J\dd{u}{x}$   then  $|\nabla u|^2$ is locally integrable
in $\Omega$, i.e $\nabla u\in L^2_{loc}(\Omega)$.

Note: In step 2, global $L^2$ estimates are obtained  on the complement
of $E$. Step 1 is used to make sure that, roughly speaking, $E$ carries
no part of the distributional gradient of $u$. This would be clear (by
abrupt cutoff and differentiation) if one knew that
$E$ has a basis of neighborhoods whose lengths of the boundaries tend to 0.
In that case Step 1 is not needed. At any rate, our approach avoids using any
non immediate property of polar sets.

\sk
{\bf Step 3.} We now want to show that for any $r>0$ (we shall need to
use $r>4$) $\nabla u\in L^r_{loc}(\Omega)$. Replacing $u$ by $u(\eta z)$
for $\eta$ small we can assume that $u$ is defined on the unit disc $\disc$
and the condition for \jhol ity will give a differential inequality
$$
  |\dd{u}{\ov z}|\leq c_0|\dd{u}{z}|, \leqno{(A)}
$$
where $c_0$ is small enough as will be said later. No more will be used.

Fix $\chi$ a cutoff function with $\chi\equiv1$ near $0$ in $\cc$ and with compact
support in
the unit disc. Let  $S$ be the convolution with the singular kernel
$\frac{-1}{\pi z^2}$.
Then we know that
$$
\dd{}{z}(\chi u)=\frac{-1}{\pi z^2}*[\dd{}{\ov z}(\chi u)]= S(\dd{}{\ov z}(\chi
u)),
$$
since $\chi u$ is in $L^2$ and of compact support.
By the Calderon-Zygmund theory for $1<r<\infty$ and $h\in L^r(\disc)$
$$
\|S(h)\|_{L^r(\disc)}=\|h\ast(\frac{-1}{\pi
z^2})\|_{L^r(\disc)}\leq\epsilon_r\|h\|_{L^r(\disc)}.
$$
For $r=2$, $\epsilon_2=1$. Write the differential inequality (A) as
$$
\dd{u}{\ov z}=-\alpha(z)\dd{u}{z},\qquad
\text{with $\alpha\in L^\infty, |\alpha|\leq c_0$}.
$$
Although $\alpha$ depends on $u$ whose regularity we want to show, consider
it to be a given function.
After cutoff
$$
\dd{\chi u}{\ov z}+\alpha(z)\dd{\chi u}{z}=g,
$$
where $g$ is a bounded function. So
$$
[{\mathbf 1}+\alpha S]\dd{\chi u}{\ov z}=g.
$$
We know that $\dd{\chi u}{\ov z}\in L^2$. So if $c_0<1$,
$[{\mathbf 1}+\alpha S]^{-1}$ is invertible on $L^2(\disc)$ and
$\dd{\chi u}{\ov z}=[{\mathbf 1}+\alpha S]^{-1}g$. But if moreover
$(|\alpha|\leq)c_0<\frac{1}{\epsilon_r}$ then this $L^2$ inverse
(given by $\sum (-1)^{k-1}(\alpha S)^k$) maps $L^r$ into $L^r$, so  $\nabla(\chi u)
\in L^r$.
(Sikorav seems to give an argument with less care.)

\sk
{\bf Step 4.} We now claim that $\Delta(\varphi\circ u)\in
L^{r/2}_{loc}(\Omega)$,
where $\varphi=\RE z_j$ or $\IM z_j$ as in Step 1.

Indeed off $E$ we have
$$
|\Delta(\varphi\circ u)|=|d\, d_J^c\varphi (\dd{u}{x},J(u)\dd{u}{x})|
\leq C|\dd{u}{x}|^2.
$$
As shown in Step 1, $\Delta(\varphi\circ u)$ is a measure on $\Omega$, we now
only need to check that it has no mass on $E$. It is given by Lemma 3.

\sk
{\bf The end.}
Fix $\infty>r>4$. Since $\Delta u\in L^{r/2}_{loc}(\Omega)$
then $\nabla u\in   C^{1-\frac{4}{r}}(\Omega)$. So $u\in C^1$ on
$\Omega$ and therefore by continuity $u$ is \jhol\ on $\Omega$.

\bigskip
\noindent
{\bf 3. Proof of Theorem B  implying  Theorem C}

\sk

Recall that we are given the following:
  $\sC$ is a closed subset of an \acm\ $(M,J)$.
Assume that there is a discrete subset   $S$ of $\sC$  such that $\sC\setminus
  S$
is locally exact $J$-pluripolar.
$\Omega$ is a connected open subset of $\cc$
and   $u\colon\Omega\to M$   is  a \contm. Assume that
$u$ is \jhol\ on $\Omega\setminus u^{-1}(\sC)$.
Assume that there exists $z_0$ with $u(z_0)\not\in\sC$. We want
to show that $u$ is \jhol.

\sk

Let $V$ be  the set of $z\in\Omega$ such that

a) $u$ is \jhol\ on a neighborhood of $z$

b) $u(z)\not\in S$ (the exceptional discrete set)

c) no neighborhood of $z$ is mapped into $\sC$.

Let  $V_0$ be the connected component of $V$ containing $z_0$.
We claim that $V_0$ is an open and relatively closed subset of the open set
$$X=\{z\in\Omega; u(z)\not\in S\}.$$
Openness requires a justification because of c). If $z_1\in V$ there exists a
neighborhood $W$ of $u(z_1)$ and a \jpshfn\ $\rho$ in $W$,
\cont\ as a map into $\{-\infty\}\cup\rr$, such that
$\rho^{-1}(-\infty)=\sC\cap W$ (possibly empty).
Then $\rho\circ u$ is  \sub\ near $z_1$ and not identical to $-\infty$.
Hence $(\rho\circ u)^{-1}(-\infty)$ has empty interior. No open subset
of some neighborhood
of $z_1$ can be mapped into $\sC$.

Now we check that $V_0$ is a relatively closed subset by the same
argument. If $z_1\in X\cap\ov V_0$, take $\rho$ as before $\rho\circ u$  is
\sub\ (because it is \sub\ when it is not $-\infty$), so the set
of $z$ such that $u(z)\in\sC$ is a   closed
polar set in a neighborhood of $z_1$. By
Theorem~B, $u$ is \jhol\ at $z_1$. So a) is proved and
c)  is trivial.

Therefore, $V_0$ is a connected component  of $X$.
If $z_2\in\Omega$ is on the boundary $bV_0$ one must have $u(z_2)\in S$, by the
definition of $X$.


Fix $z_2\in \Omega\cap bV_0$.
Since $S$ is discrete, there exists $\epsilon>0$ such that $u(z)\equiv p\in S$
for all $z\in bV_0$ such that $|z-z_2|<\epsilon$. Let $L$ be a      \jpshfn\ defined
near
 $u(z_2)$ with pole $(-\infty)$ at  $u(z_2)$   (Chirka's function, \cont\ away from $u(z_2)$).
Shrinking
$\epsilon$ if needed we can assume that
 $r(z)=L\circ u(z)$
is defined on the disc $\{z\in\cc; |z-z_2|<\epsilon\}$, where on that disc
$$
r(z)=L\circ u(z)\ \text{if $z\in V_0$}, \quad r(z)=-\infty\ \text{if $z\not\in
V_0$}.
$$
$r$ is a \subfn\ and $r^{-1}(-\infty)$ is a   closed
polar subset, which is removable  for
$u$   by
Theorem~B.

 Therefore, $u$ is $J$-holomorphic at each point of $\Omega\cap bV_0$
(and $J$-holomorphic on $V_0$ by definition).
And
$\Omega\setminus V_0$ has empty interior. Otherwise, by the connectedness of $\Omega$
 the interior of $\Omega\setminus V_0$ has a boundary point
  $z_2\in bV_0\cap\Omega$,  and  the above polar set
 $r^{-1}(-\infty)$ has non-empty interior, which is a contradiction.
 This shows that $u$ is $J$-holomorphic on $\Omega$.

\newcommand{\BCfise}{\bibitem[B-C]{BCfise} F. Bruhat and H. Cartan,
{\em Sur la structure des sous-ensembles analytiques r\'eels},
C. R. Acad. Sci. Paris {\bf 244}(1957), 988--991.
}

\newcommand{\ISnini}{\bibitem[I-S]{ISnini}
S. Ivashkovich and V. Shevchishin, {\em Structure of the moduli space
in a neighborhood of a cusp-curve and meromorphic hulls}, Invent. Math.  {\bf
136}(1999),  no. 3, 571--602.
}

\newcommand{\Pazefi}{\bibitem[P]{Pazefi}
N. Pali, {\em Faisceaux $\db$-coh\'erents sur les vari\'et\'es complexes},
preprint.
}

\newcommand{\Ahsisi}{\bibitem[A]{Ahsisi}
L.V. Ahlfors, {\em Lectures on quasiconformal mappings},
  Van Nostrand Mathematical Studies, No. 10 D. Van Nostrand Co., Inc., Toronto,
Ont.-New York-London 1966.
}

\newcommand{\Sinifo}{\bibitem[S]{Sinifo}
  J.-C. Sikorav, {\em Some properties of holomorphic curves in almost complex
  manifolds}, in
  Holomorphic curves in symplectic geometry,  165--189, Progr. Math., 117,
Birkh\"auser, Basel, 1994.
  }

\newcommand{\Deeitw}{\bibitem[D]{Deeitw}
J.-P. Demailly, {\em
Estimations $L\sp{2}$ pour l'op\'erateur $\bar \partial $ d'un fibr\'e
vectoriel holomorphe
semi-positif au-dessus d'une vari\'et\'e K\"ahl\'erienne compl\`ete},
  Ann. Sci. \'Ecole Norm. Sup. (4)  {\bf 15}(1982),  no. 3, 457--511.}

\newcommand{\FHSnifi}{\bibitem[F-H-S]{FHSnifi} A. Floer,  H. Hofer, and D.
Salamon, {\em
  Transversality in elliptic Morse theory for the symplectic action},
   Duke Math. J.  {\bf 80}(1995),  no. 1, 251--292.}

   \newcommand{\IRzefo}{\bibitem[I-R]{IRzefo}
S. Ivashkovich and J.-P. Rosay, {\em
  Schwarz-type lemmas for solutions of
  $\overline\partial$-inequalities and complete
  hyperbolicity of almost complex manifolds},
  Ann. Inst. Fourier (Grenoble)  54  (2004),  no. 7, 2387--2435.}

  \newcommand{\Rozesi}{\bibitem[R]{Rozesi}
  J.-P. Rosay, {\em  J-holomorphic submanifolds are pluripolar}, to appear
  in Math. Zeit.}

   \newcommand{\Nasisi}{\bibitem[N]{Nasisi}
  R.  Narasimhan, {\em Introduction to the theory of analytic spaces},
  Lecture Notes in Mathematics, No. 25 Springer-Verlag, Berlin-New York 1966.
  }

\newcommand{\DFseei}{\bibitem[D-F]{DFseei} K. Diederich and J.E. Forn{\ae}ss,
{\em Pseudoconvex domains with real-analytic boundary},
Ann. Math. (2) {\bf 107}(1978), no. 2, 371--384.
}

  \newcommand{\MDSnifo}{\bibitem[MD-S]{MDSnifo}
D. McDuff and D. Salamon, {\em $J$-holomorphic curves and
quantum cohomology}, University Lecture Series, 6. American
Mathematical Society, Providence, RI, 1994. }

\newcommand{\HLeifo}{\bibitem[H-L]{HLeifo}
 G. Henkin and J. Leiterer, {\em Theory of functions on complex manifolds},
 Monographs in Mathematics, 79. Birkh\"{a}user Verlag, Basel, 1984.}

\end{document}